\documentclass{amsart}

\usepackage{amsfonts}
\usepackage{cases}
\usepackage{mathrsfs}
\usepackage{bbm}
\usepackage{amssymb}
\usepackage{txfonts}
\usepackage{amscd}
\usepackage{amsfonts,latexsym,amsmath,amsthm,amsxtra,mathdots}
\usepackage[bookmarks,colorlinks]{hyperref}
\usepackage[all,cmtip]{xy}
\RequirePackage{amsmath} \RequirePackage{amssymb}
\usepackage{color}
\usepackage{colordvi}
\usepackage{multicol}
\usepackage{pdfsync}
\usepackage[utf8]{inputenc}
\usepackage{hyperref}
\usepackage{graphicx}
\usepackage{amsmath}
\usepackage{amsmath,amscd}
\usepackage[normalem]{ulem}
\usepackage{textcomp}
\usepackage{enumitem}
\usepackage{tikz}
\usetikzlibrary{matrix,arrows,trees}
\usetikzlibrary{positioning}
\usepackage{pdfcolmk}
\usepackage{todonotes}

\newcommand{\cells}{\diamond\,}

\newcommand{\id}{{{\rm id}}}

    \newcommand{\BA}{{\mathbb {A}}}

     \newcommand{\BJ}{{\mathbb {J}}}
     
     \newcommand{\BN}{{\mathbb {N}}}
     
    \newcommand{\BQ}{{\mathbb {Q}}}

     \newcommand{\BZ}{{\mathbb {Z}}}

     \newcommand{\CF}{{\mathcal {F}}}

\def\-{^{-1}}

\newcommand{\delete}[1]{}

    \theoremstyle{plain}

\newtheorem{thm}{Theorem}[section]
\newtheorem{defn}[thm]{Definition} 
\newtheorem{ex}[thm]{Example} 
\newtheorem{lem}[thm]{Lemma}
\newtheorem{prop}[thm]{Proposition}

\newtheorem{rem}[thm]{Remark}

    \numberwithin{equation}{section}

\def\Proof{\noindent{\bf Proof}\quad}
\def\qed{\hfill$\square$\smallskip}

\newcommand{\abs}[1]{\lvert#1\rvert}

\begin{document}

\title{The  Isomorphism Conjecture for solvable groups in  Waldhausen's A-theory}

\author{F. Thomas Farrell}
\address{Department of Mathematical Sciences and Yau Mathematical Sciences
Center, Tsinghua University, Beijing, China}
\email{farrell@math.binghamton.edu}

\author{Xiaolei Wu}
\address{Max Planck Institute for Mathematics, Vivatsgasse 7, 53111 Bonn, Germany }
\email{hsiaolei.wu@mpim-bonn.mpg.de}



\date{September, 2017}


\begin{abstract}

We prove the A-theoretic Isomorphism Conjecture with coefficients and finite wreath products  for solvable groups.
\end{abstract}

\maketitle

\section{introduction}

Recently, there has been a growing interest in applying techniques used to prove the Isomorphism Conjecture in K- and L- theory (see for example \cite{Bartels(2012),BFL,BL1, BL12,BLR2,BLRH}) to the A-theory setting. In particular, Ullman--Winges generalized Farrell-Hsiang method \cite{BL1,FH} to the A-theory setting and used it to verify the A-theory Isomorphism Conjecture for virtually poly-$\BZ$-groups \cite{UW}. Later, together with Enkelmann-- L\"uck--Pieper, they proved if a group is homotopy transfer reducible, then it satisfies the A-theory Isomorphism Conjecture \cite{ELPUW}. In particular, they proved the conjecture for hyperbolic and CAT(0) groups.

Following \cite[Conjecture 7.1]{UW}(see also \cite[Conjecture 2.12]{ELPUW}), we formulate the Isomorphism Conjecture in A-theory as follows.  Let $G$ be a discrete group and $W$ be a connected $G$-CW-complex with free  $G$-action. Let $\BA^{-\infty}(W)$ be the non-connective delooping of Waldhausen's algebraic $K$-theory of $W$ over the orbit category. The \emph{Isomorphism Conjecture} predicts that the \emph{assembly map}
\begin{equation}
 H_n^G(E_{\mathcal{VC} yc} G; \BA^{-\infty}_W) \to H_n^G(G/G;\BA^{-\infty}_W) \cong \pi_n\BA^{-\infty}(G\backslash W),
\end{equation}
which is induced by the projection map $E_{\mathcal{VC} yc} G \to G/G$ from the classifying space for virtually cyclic subgroups to a point, is an isomorphism for all $n \in \BZ$. We say the A-theory isomorphism Conjecture with coefficients holds for $G$ if it holds for any  such $W$. We say the A-theory Isomorphism Conjecture with coefficients and finite wreath products holds for $G$ if the conjecture with coefficients holds for any  wreath products $G \wr F$, where $F$ is a finite group. Our main theorem can now be stated as follows.
\begin{thm}\label{main theorem}
The Isomorphism Conjecture in A-theory with coefficients and finite wreath products holds for all solvable groups.
\end{thm}
\begin{rem}
As pointed out in \cite[Remark 1.4]{ELPUW}, our result combined with \cite[Theorem 1.1 + Corollary 6.6]{ELPUW} and arguments in  \cite{KLP,Rup} implies that the Isomorphism Conjecture in A-theory with coefficients and finite wreath products  holds for any (not necessarily cocompact)
  lattice in a second countable locally compact Hausdorff group with finitely many path
  components, the groups $GL_n(\BQ)$ and $GL_n(F(t))$ for $F(t)$ the function field over a
  finite field $F$, and all $S$-arithmetic groups.
\end{rem}

\begin{rem}
Theorem \ref{main theorem} combined with \cite[Theorem 1.1]{ELPUW} and arguments due to Gandini-- Meinert--R\"uping in \cite{GMR}   imply that any fundamental groups of graphs of abelian groups satisfy the Isomorphism Conjecture in A-theory with coefficients and finite wreath products. In particular  all Baumslag-Solitar groups satisfy the conjecture.
\end{rem}

\begin{rem}
Independently,  Kasprowski--Ullmann--Wegner--Winges \cite{DUW} also obtained a proof for Theorem  \ref{main theorem}.
\end{rem}

Our proof is based on Wegner's paper \cite{W2}. The proof can not be directly taken over due to the fact that one has to change hyper-elementary subgroups in Farrell--Hsiang method (for K- and L- theory \cite{BL1}) to Dress subgroups (for A-theory \cite{UW}).  The paper is organized as follows. In Section \ref{section-DFHJ}, we generalize Wegner's Farrell-Hsiang-Jones group to the A-theory setting and prove it satisfies the A-theory Isomorphism Conjecture. In section \ref{section-dress}, we embed the group $\BZ[w,\frac{1}{w}]\rtimes_w \BZ$\footnote{multiplication in $\BZ[w,\frac{1}{w}]\rtimes_w \BZ$ is given by  $(x_1,y_1)(x_2,y_2) = (x_1 + w^{y_1}x_2,y_1+y_2)$,
where $(x_i,y_i) \in\BZ[w,\frac{1}{w}] \rtimes_w \BZ$.}  (we only need to study these groups due to Lemma \ref{reduction}) into a slightly bigger group and study finite quotients of it. We show that the Dress subgroups inside these finite quotients have large index. In the last section, we prove our main theorem based on work of Wegner \cite{W2}.

In this paper, when we say a group satisfies the full Isomorphism Conjecture in A-theory we mean it satisfies the Isomorphism Conjecture with coefficients and finite wreath products in A-theory.

\textbf{Acknowledgements.}  The second author is supported by the Max Planck Institute for Mathematics at Bonn.  We want to thank  Guoliang Yu for inviting us to visit the  Shanghai Center for Mathematical Sciences, where this project was initiated. The second author also want to thank  Yang Su for providing accommodation support during his visit in Beijing and Malte Pieper for very helpful discussions on \cite[Section 6-7]{ELPUW}.

\section{Dress--Farrell--Hsiang--Jones group} \label{section-DFHJ}

In this section, we generalize Wegner's Farrell--Hsiang--Jones group to Dress-Farrell--Hsiang--Jones group and show that it  satisfies the Isomorphism Conjecture in  A-theory. Dress--Farrell--Hsiang--Jones group  is a combination of the Dress--Farrell--Hsiang group  \cite{UW} and homotopy transfer reducible group \cite[Defintion 6.2]{ELPUW}.

\subsection{Dress-Farrell-Hsiang group}

We follow Ullmann and Winges' terminology \cite{UW}, define Dress group as follows,
\begin{defn}\label{dress}
A finite group $D$ is  a Dress
group if there are primes $p_1$ and $p_2$ and a normal series $D_1 \lhd D_2 \lhd D$ such that $D_1$ is a
$p_1$-group, $D_2/D_1$ is cyclic and $D/D_2$ is a $p_2$-group.
\end{defn}
Note that the primes $p_1$ and $p_2$ in Definition \ref{dress} need not be distinct
\begin{lem}\label{dressdef}
In the definition of Dress group, we can always assume the order  satisfies
$$(|D_1|, |D_2/D_1|) = (|D/D_2|,|D_2/D_1|) = 1.$$
In particular $D_2 \cong D_1 \rtimes \BZ/|D_2/D_1|$.
\end{lem}

\Proof We first show we can assume $(|D_1|, |D_2/D_1|) =1$. Suppose $(|D_1|, |D_2/D_1|) \neq 1$, then it must be a power of $p_1$ since $D_1$ is a $p_1$-group. Assume $(|D_1|, |D_2/D_1|) = p_1^k$ for some $k \geq 1$ and let the projection from $D_2$ to $D_2/D_1$ by $\pi_1$. Then $D_1'=\pi_1^{-1}(\frac{|D_2/D_1|}{p_1^k} D_2/D_1)$ is a normal subgroup of $D_2$ such that $(|D_1'|, |D_2/D_1'|) = 1$, $D_1'$ is a
$p_1$-group, $D_2/D_1'$ is cyclic. Since $D_2/D_1'$ is cyclic and $(|D_1'|, |D_2/D_1'|) = 1$, the projection map $D_2\rightarrow D_2/D_1'$ splits and we have $D_2 \cong D_1' \rtimes \BZ/|D_2/D_1'|$.

Now we assume $D_2 \cong D_1 \rtimes \BZ/|D_2/D_1|$ and $(|D_1|, |D_2/D_1|) = 1$. Note that $D_1$ is a characteristic subgroup of $D_2$. Hence it is a normal subgroup of $D$. We have the following short exact sequence:

$$\begin{CD}
1    @>  >> D_2/D_1 @>>> D/D_1  @>>> D/D_2 @>>> 1\\
\end{CD}$$

Now $D_2/D_1$ is a cyclic group and $D/D_2$ is a $p_2$-group. If $(|D/D_2|,|D_2/D_1|) \neq 1$, then, $(|D/D_2|,|D_2/D_1|) = p_2^l$ for some $l>0$. This means $|D_2/D_1| = p_2^l s$ where $(p_2,s)=1$. Now  $|D_2/D_1| = \BZ/( p_2^l s) \cong \BZ/p_2^l \times \BZ/s$. Since $\BZ/p_2^l$ and $\BZ/s$ are characteristic subgroup of $\BZ/( p_2^l s)$, we have a characteristic subgroup of order $s$ in $D_2/D_1$ which is normal in $D/D_2$. Let the corresponding group in $D_2$ be $D_2'$, then we would have $(|D/D_2'|,|D_2'/D_1|) = 1$. The new $D_2'$ is a subgroup of $D_2$ which also has the property that $(D_1',D_2'/D_1') =1$ for the corresponding $D_1'$ constructed in the first paragraph.

\qed

Recall the definition of the \emph{$\ell^1$-metric} on a simplicial complex. If $X$ is a simplicial complex and $\xi = \sum_x \xi_x \cdot x$, $\eta = \sum_x \eta_x \cdot x$ are points in $X$, this metric is given by $d^{\ell^1}(\xi,\eta) = \sum_x \abs{\xi_x - \eta_x}$. All simplicial complexes we consider are equipped with this metric.

We call a generating set $S$ of a group $G$ \emph{symmetric} if $s \in S$ implies $s^{-1} \in S$. The following group is defined by Ullmann and Winges in \cite[Section 1]{UW}.

\begin{defn}\label{def_dfhgroup}
 Let $G$ be a group and $S$ be a symmetric, finite generating set of $G$. Let $\mathcal{F}$ be a family of subgroups of $G$.

 Call $(G,S)$ a \emph{Dress--Farrell--Hsiang group with respect to $\mathcal{F}$} if there exists $N \in \BN$
 such that for every $\epsilon > 0$ there is an epimorphism $\pi \colon G \twoheadrightarrow F$ onto a finite group $F$
 such that the following holds: For every Dress subgroup $D \subset F$, there are
 a $\overline{D} := \pi^{-1}(D)$-simplicial complex $E_D$ of dimension at most $N$ whose isotropy groups lie in $\mathcal{F}$, and
 a $\overline{D}$-equivariant map $\phi_D \colon G \to E_D$ such that $d^{\ell^1}(\phi_D(g),\phi_D(g')) \leq \epsilon$
 whenever $g^{-1}g' \in S$.
\end{defn}

Ullmann and Winges proves the following \cite[Theorem 7.4]{UW}.

\begin{thm}\label{thm_dfhmethod}
 Let $G$ be a discrete group. Suppose that there is a symmetric, finite generating set
 $S \subset G$ and a family of subgroups $\CF$ of $G$ such that $(G,S)$ is a Dress--Farrell--Hsiang group with respect to $\CF$.
 Then $G$ satisfies the Isomorphism Conjecture in $A$-theory with coefficients  with respect to $\CF$.
\end{thm}

\subsection{Homotopy transfer reducible groups}

 \begin{defn}\label{def:hptycohG-action} \cite[Definition 6.1]{ELPUW}
  A \emph{homotopy coherent $G$-action} of a group $G$ on a topological space $X$ is a continuous map
  \[
  \Psi: \coprod_{j=0}^\infty((G \times [0,1])^j\times G\times X)\to X
  \]
  with the following properties:
  \begin{equation*}
       \Psi(\gamma_k,t_k,\dots,\gamma_1,t_1,\gamma_0,x) = \begin{cases}
       \Psi(\dots, \gamma_j,  \Psi(\gamma_{j-1},\dots,x)) & t_j = 0 \\
        \Psi(\dots,\gamma_j\gamma_{j-1},\dots,x) & t_j = 1 \\
        \Psi(\gamma_k,\dots,\gamma_2,t_2,\gamma_1,x) & \gamma_0 = e \\
        \Psi(\gamma_k,\dots,t_{j+1}t_j,\dots,\gamma_0,x) & \gamma_j = e, 1 \leq j < k \\
        \Psi(\gamma_{k-1},t_{k-1},\dots,t_1,\gamma_0,x) & \gamma_k = e\\
        x & \gamma_0 =e, k=0
      \end{cases}
  \end{equation*}
  \end{defn}

   \begin{defn}\cite[Definition 6.8]{ELPUW}\label{def-metric}
   Let $(X,d_X)$ be a metric space, $ \Psi$ a homotopy coherent $G$-action on $X$, and $S \subset G$ a finite subset containing the trivial element. Let $k \in \BN$ and $\Lambda > 0$. Define
    \begin{equation*}
     d_{S,k,\Lambda}( (x,g), (y,h) ) \in [0,\infty]
    \end{equation*}
   to be the infimum over the numbers
    \begin{equation*}
     l + \sum_i \Lambda \cdot d_X(x_i,z_i),
    \end{equation*}
   where the infimum is taken over all $l \in \BN$, $x_0,\dots,x_l$, $z_0,\dots,z_l \in X$ and $a_1,\dots,a_l$, $b_1,\dots,b_l \in S$ such that
    \begin{enumerate}
     \item $x_0 = x$ and $z_l = y$;
     \item $ga_1^{-1}b_1\dots a_l^{-1}b_l = h$;
     \item for each $1 \leq i \leq l$ there are elements $r_0,\dots,r_k,s_0,\dots,s_k \in S$ such that $a_i = r_k\dots r_1$, $b_i = s_k \dots s_1$ and $ \Psi(r_k,t_k,\dots,r_0,z_{i-1}) =  \Psi(s_k,u_k,\dots,s_0,x_i)$ for some $t_1, \dots,t_k,u_1,\dots,u_k \in [0,1]$.
    \end{enumerate}
If no such data exist, take the infimum to be $\infty$.
   \end{defn}

  \begin{defn}\label{def:transfer-reducible}
   Let $G$ be a discrete group. Let $S \subset G$ be a finite, symmetric generating set of $G$ which contains the trivial element and $S^n$ be the ball of radius $n$ around the identity element under the word metric of $G$ with respect to $S$. Let $\CF$ be a family of subgroups of $G$.

   Then $G$ is \emph{homotopy transfer reducible over $\CF$} if there exists $N \in \BN$ such that  for every $n \in \BN$ there are
    \begin{enumerate} [label=(\roman*)]
     \item a compact, contractible metric space $(X,d_X)$ such that for every $\epsilon > 0$ there is an $\epsilon$-controlled domination of $X$ by an at most $N$-dimensional, finite simplicial complex.
     \item a homotopy coherent $G$-action $\Psi$ on $X$.
     \item a $G$-simplicial complex $\Sigma$ of dimension at most $N$ whose isotropy is contained in $\CF$.
      \item a positive real number $\Lambda$;
     \item a $G$-equivariant map $\phi \colon X \times G \to \Sigma$ such that
     \[
     n \cdot d^{\ell^1}(\phi(x,g),\phi(y,h)) \leq d_{S^n,n,\Lambda}((x,g),(y,h))
     \]
     holds for all $(x,g)$, $(y,h) \in X \times G$.
    \end{enumerate}
  \end{defn}

\begin{rem}\label{diff-transfer-reducible}
Note our definition of homotopy transfer reducible differs slightly from the one in \cite[Definition 6.2]{ELPUW}, but as in the proof of \cite[Theorem 6.19]{ELPUW}  shows, our definition implies the definition there.
\end{rem}

  \begin{thm}\label{thm:afjc-transfer-reducible}\cite[Theorem 6.14]{ELPUW}
   Let $G$ be a discrete group and let $\CF$ be a family of subgroups of $G$.
   If $G$ is homotopy transfer reducible over $\CF$, then $G$ satisfies the Isomorphism Conjecture with coefficients in A-theory with respect to $\CF$.
  \end{thm}

\subsection{Category of controlled retractive $G$-CW-complexes}
  Let $G$ be a discrete group and  $\CF$ be a family of subgroups of $G$.
  In \cite[Corollary 6.11]{UW} it was shown that the Isomorphism Conjecture with coefficients holds for $G$ iff the algebraic K-theory of $\mathfrak{R}^G_f(W,\BJ(G,E_{\mathcal{F}}(G)),h)$  vanishes for every free $G$-CW-complex $W$, so it plays the role of obstruction category as in the proof of the Isomorphism Conjecture for algebraic K-theory, see for example \cite[Section~3]{BLR2}.  Most of the material here are directly borrowed from \cite[Section 6.2]{ELPUW} which in turn come from  \cite{UW}.

  \begin{defn}\label{coarse-structure}\cite[Definition 2.1]{UW}
 Let $Z$ be a $G$-space which is Hausdorff.
 A set of \emph{morphism control conditions} $\mathfrak{Z}$ is a collection of $G$-invariant subsets of $Z \times Z$
 with the following properties:
 \begin{enumerate}[label={(C{\arabic*})}]
  \item  Every $C \in \mathfrak{C}$ contains the diagonal $\Delta(Z) := \{ (z,z) \mid z \in Z \}$.
  \item  Every $C \in \mathfrak{C} $ is symmetric.
  \item  For all $C, C' \in \mathfrak{C}$ there is some $C'' \in \mathfrak{C}$ such that $C \cup C' \subset C''$.
  \item  For all $C, C' \in \mathfrak{C}$ there is some $C'' \in \mathfrak{C}$ such that $C' \circ C \subset C''$,
   where the composition $C' \circ C$ is defined as
   \begin{equation*}
    C' \circ C := \{ (z'',z) \mid \exists z' \colon (z',z) \in C, (z'',z') \in C' \}.
   \end{equation*}
 \end{enumerate}
 A set of \emph{object support conditions} $\mathfrak{S}$ is a collection of $G$-invariant subsets of $Z$
 with the following properties:
 \begin{enumerate}[label={(S\arabic*})]
   \item  For all $S, S' \in \mathfrak{S}$ there is some $S'' \in \mathfrak{S}$ such that $S \cup S' \subset S''$.
 \end{enumerate}
 The triple $\mathfrak{Z} = (Z,\mathfrak{C},\mathfrak{S})$ is called a \emph{coarse structure}.
\end{defn}

  Fix a coarse structure $\mathfrak{Z}$.  A \emph{labeled $G$-CW-complex relative W} \cite[Definition~2.3]{UW}, is a pair $(Y, \kappa)$, where $Y$ is a free $G$-CW-complex relative $W$ together with a $G$-equivariant function $\kappa \colon \cells Y \rightarrow Z$. Here, $\diamond Y$ denotes the (discrete) set of relative cells of $Y$.

  A \emph{$\mathfrak{Z}$-controlled map} $f \colon (Y_1,\kappa_1) \rightarrow (Y_2, \kappa_2)$ is a $G$-equivariant, cellular map $f \colon Y_1 \rightarrow Y_2$ relative $W$ such that for all $k \in \mathbb{N}$ there is some $C \in \mathfrak{C}$ for which
  \[
  (\kappa_2,\kappa_1)(\{(e_2,e_1) \mid e_1 \in \cells_k Y_1, e_2\in \cells Y_2, \langle f(e_1) \rangle \cap e_2 \neq \emptyset\}) \subseteq C
  \]
  holds, where $\cells_k Y_1$ denotes the set of relative $k$-cells of $Y$.

  A \emph{$\mathfrak{Z}$-controlled $G$-CW-complex relative W} is a labeled $G$-CW-complex $(Y,\kappa)$ relative $W$, such that the identity is a $\mathfrak{Z}$-controlled map and for all $k \in \BN$ there is some $S \in \mathfrak{S}$ such that
  \[
    \kappa(\diamond_k Y)\subseteq S.
  \]

  A \emph{$\mathfrak{Z}$-controlled retractive space relative $W$} is a $\mathfrak{Z}$-controlled $G$-CW-complex  $(Y,\kappa)$ relative $W$ together with a $G$-equivariant retraction $r \colon Y \to W$, i.e., a left inverse to the  inclusion $W \hookrightarrow Y$. The $\mathfrak{Z}$-controlled retractive spaces relative $W$ form a category $\mathcal{R}^G(W,\mathfrak{Z})$ in which morphisms are $\mathfrak{Z}$-controlled maps which additionally respect the chosen retractions.

  The category of controlled $G$-CW-complexes (relative $W$) and controlled maps admits a notion of \textit{controlled homotopies}\cite[Definition~2.5]{UW} via the objects $(Y \leftthreetimes [0,1], \kappa \circ pr_Y)$, where $Y \leftthreetimes [0,1]$ denotes the reduced product which identifies $W \times [0,1] \subseteq Y \times [0,1]$ to a single copy of $W$ and $pr_Y: \cells Y \leftthreetimes [0,1] \rightarrow \diamond Y$ is the canonical projection.
  In particular, we obtain a notion of \emph{controlled homotopy equivalence} (or \emph{$h$-equivalence}).

  A $\mathfrak{Z}$-controlled retractive space $(Y, \kappa)$ is called \textit{finite} \cite[Definition~3.3]{UW} if it is finite-dimensional, the image of $Y \backslash W$ under the retraction meets the orbits of only finitely many path components of $W$ and for all $z \in Z$ there is some open neighborhood $U$ of $z$ such that $\kappa^{-1}(U)$ is finite.

  A $\mathfrak{Z}$-controlled retractive space $(Y, \kappa)$ is called \textit{finitely dominated}, if there are a finite $\mathfrak{Z}$-controlled, retractive $G$-CW-complex $D$ relative to $W$, a morphism $p \colon D \rightarrow Y$ and a $\mathfrak{Z}$--controlled map $i \colon Y \rightarrow D$ such that $p \circ i \simeq_{\mathfrak{Z}} \id_Y$ as $\mathfrak{Z}$-controlled maps.

  The finite, respectively finitely dominated, $\mathfrak{Z}$-controlled retractive spaces form full subcategories $\mathcal{R}^G_f(W,\mathfrak{Z}) \subset \mathcal{R}^G_{fd}(W,\mathfrak{Z}) \subset \mathcal{R}^G(W,\mathfrak{Z})$.
  All three of these categories support a Waldhausen category structure in which inclusions of $G$-invariant subcomplexes up to isomorphism are the cofibrations and controlled homotopy equivalences are the weak equivalences \cite[Corollary~3.22]{UW}.

 Let $\mathfrak{Z}_1 = (Z_1, \mathfrak{C}_1, \mathfrak{S}_1)$, $\mathfrak{Z}_2 = (Z_2, \mathfrak{C}_2, \mathfrak{S}_2)$ be two coarse structures.
 A \emph{morphism of coarse structures} \cite[Definition 3.23]{UW} $z \colon \mathfrak{Z}_1 \to \mathfrak{Z}_2$ is a $G$-equivariant map of sets
 $z \colon Z_1 \to Z_2$ satisfying the following properties:
 \begin{enumerate}
  \item For every $S_1 \in \mathfrak{S}_1$, there is some $S_2 \in \mathfrak{S}_2$ such that $z(S_1) \subset S_2$.
  \item  For every $S \in \mathfrak{S}_1$ and $C_1 \in \mathfrak{C}_1$, there is some $C_2 \in \mathfrak{C}_2$ such that $(z \times z)((S \times S) \cap C_1) \subset C_2$.
  \item For every $S \in \mathfrak{S}_1$ and all subsets $A \subset S$ which are locally finite in $Z_1$, the set $z(A)$ is locally finite in $Z_2$ and for all $x \in z(A)$, the set $z^{-1}(x) \cap A$ is finite.
 \end{enumerate}

Morphism of coarse structures induce morphism of corresponding categories of controlled $G$-CW-complexes (relative to W), see \cite[Proposition 3.24]{UW} for more details.

  Let $X$ be a $G$-CW-complex and let $M$ be a metric space with free, isometric $G$-action. Define $\mathfrak{C}_{bdd}(M)$ to be the collection of all subsets $C \subset M \times M$ which are of the form
  \begin{equation*}
      C = \{ (m,m') \in M \times M \mid d(m,m') \leq \alpha \}
  \end{equation*}
  for some $\alpha \geq 0$. Define further $\mathfrak{C}_{Gcc}(X)$ to be the collection of all $C \subset (X \times [1,\infty[) \times (X \times [1,\infty[)$ which satisfy the following:
  \begin{enumerate}
      \item For every $x \in X$ and every $G_x$-invariant open neighborhood $U$ of $(x,\infty)$ in $X \times [1,\infty]$, there exists a $G_x$--invariant open neighborhood $V \subset U$ of $(x,\infty)$ such that $(((X \times [1,\infty[) \setminus U) \times V) \cap C = \emptyset$.
      \item Let $p_{[1,\infty[} \colon X \times [1,\infty[ \to [1,\infty[$ be the projection map. Equip $[1,\infty[$ with the Euclidean metric. Then there exists some $B \in \mathfrak{C}_{bdd}([1,\infty[)$ such that $C \subset p^{-1}_{[1,\infty[}(B)$.
      \item $C$ is symmetric, $G$--invariant and contains the diagonal.
  \end{enumerate}
  Next define $\mathfrak{C}(M,X)$: Let $p_M \colon M \times X \times [1,\infty[ \to M$ and $p_{X \times [1,\infty[} \colon M \times X \times [1,\infty[ \to X \times [1,\infty[$ denote the projection maps. Then $\mathfrak{C}(M,X)$ is the collection of all subsets $C \subset (M \times X \times [1,\infty[)^2$ which are of the form
  \begin{equation*}
   C = p_M^{-1}(B) \cap p_{X \times [1,\infty[}^{-1}(C')
  \end{equation*}
  for some $B \in \mathfrak{C}_{bdd}(M)$ and $C' \in \mathfrak{C}_{Gcc}(X)$.

  Finally, define $\mathfrak{S}(M,X)$ to be the collection of all subsets $S \subset M \times X \times [1,\infty[$ which are of the form $S = K \times [1,\infty[$ for some $G$-compact subset $K \subset M \times X$.

  All these data combine to a coarse structure
  \begin{equation*}
      \BJ(M,X) := (M \times X \times [1,\infty[, \mathfrak{C}(M,X), \mathfrak{S}(M,X))
  \end{equation*}
  which serves to define the ``obstruction category" $\mathfrak{R}^G_f(W,\BJ(G,E_{\mathcal{F}}(G)),h)$~\cite[Example~2.2 and Definition~6.1]{UW}, where $E_{\mathcal{F}}(G)$ is the classifying space for $G$ respect to the family of subgroups $\mathcal{F}$. We will consider the non-connective $K$-theory spectrum of $\mathfrak{R}^G_f(W,\BJ(G,E_{\mathcal{F}}(G)))$ with respect to the $h$-equivalences\cite[Section~5]{UW}.

  \begin{prop}\label{obstruction}\cite[Corollary~6.11]{UW} A group $G$ satisfies the Isomorphism  Conjecture with coefficients in $A$-theory with respect to $\mathcal{F}$ if and only if the algebraic K-theory of $\mathfrak{R}^G_f(W,\BJ(G,E_{\mathcal{F}}(G)),h)$ vanishes for every free $G$-CW-complex $W$.
\end{prop}

Suppose that $(M_n)_n$ is a sequence of metric spaces with a free, isometric $G$-action. Let $X$ be a $G$-CW-complex. Following \cite[Section~7]{UW}, define the coarse structure
  \begin{equation*}
      \BJ((M_n)_n,X) := \big( \coprod_n M_n \times X \times [1,\infty[, \mathfrak{C}((M_n)_n,X), \mathfrak{S}((M_n)_n,X) \big)
  \end{equation*}
  as follows: Members of $\mathfrak{C}((M_n)_n,X)$ are of the form $C = \coprod_n C_n$ with $C_n \in\mathfrak{C}(M_n,X)$, and we additionally require that $C$ satisfies the \emph{uniform metric control condition}: There is some $\alpha > 0$, independent of $n$, such that for all $((m,x,t)$, $(m',x',t')) \in C_n$ we have $d(m,m') < \alpha$. Members of $\mathfrak{S}((M_n)_n,X)$ are sets of the form $S = \coprod_n S_n$ with $S_n \in \mathfrak{S}(M_n,X)$. The resulting category $\mathfrak{R}^G(W,\BJ((M_n)_n,X))$ is canonically a subcategory of the product category $\prod_n \mathfrak{R}^G(W,\BJ(M_n,X))$.

We need also another class of weak equivalences on $\mathfrak{R}^G(W,\BJ((M_n)_n,E))$ to describe the target of the transfer in the proof of Theorem \ref{DFHJ-iso}. These $h^{fin}$-equivalences were introduced in the proof of \cite[Theorem 10.1]{UW}. Let $(M_n)_n$ be a sequence of metric spaces with free, isometric $G$-action. Let $(Y_n)_n$ be an object of $\mathfrak{R}^G(W, \BJ((M_n)_n,E))$. For $\nu \in \BN$, we denote by $(-)_{n > \nu}$ the endofunctor which sends $(Y_n)_n$ to the sequence $(\widetilde{Y}_n)_n$ with $\widetilde{Y}_n = \ast$ for $n \leq \nu$ and $\widetilde{Y}_n = Y_n$ for $n > \nu$. A morphism $(f_n)_n \colon (Y_n)_n \rightarrow (Y_n')_n$ is an \emph{$h^{fin}$-equivalence} if there is some $\nu \in \BN$, such that $(f_n)_{n>\nu}\colon (Y_n)_{n>\nu} \rightarrow (Y_n')_{n>\nu}$ is an $h$-equivalence. To distinguish them, we will denote the one with $h$-equivalence by $\mathfrak{R}^G(W,\BJ((M_n)_n,E),h)$ and the new one with   $h^{fin}$-equivalence by $\mathfrak{R}^G(W,\BJ((M_n)_n,E),h^{fin})$.

\subsection{Dress--Farrell--Hsiang--Jones group}
In this subsection we introduce and study Dress-Farrell-Hsiang-Jones groups. They arise from a combination of the Farrell-Hsiang method and transfer reducibility.

\begin{defn}[DFHJ group] \label{def-DFHJ}
Let $G$ be a finitely generated group and let $\mathcal{F}$ be a family of subgroups. Let $S \subseteq G$ be a finite symmetric subset which generates $G$ and contains the trivial element $e \in G$.
We call $G$ a \emph{Dress--Farrell--Hsiang--Jones  group  with respect to the family $\mathcal{F}$} if there exist a natural number $N$ and surjective homomorphisms $\alpha_n \colon G \to F_n$ (for any $n \in \BN$) onto finite groups $F_n$ such that the following condition is satisfied.
For any Dress subgroup $D$ of $F_n$ there exist
\begin{itemize}
\item a compact, contractible, controlled $N$-dominated metric space $X_{n,D}$,
\item a  homotopy coherent $G$-action $\Psi_{n,D}$ on $X_{n,D}$,
\item a positive real number $\Lambda_{n,D}$,
\item a simplicial complex $E_{n,D}$ of dimension at most $N$ with a simplicial $\alpha_n^{-1}(D)$-action whose stabilizers belong to $\mathcal{F}$,
\item and an $\alpha_n^{-1}(D)$-equivariant map $f_{n,D} \colon G \times X_{n,D} \to E_{n,D}$
\end{itemize}
such that
\[
n \cdot d_{E_{n,D}}^1\big(f_{n,D}(g,x),f_{n,D}(h,y)\big) \leq d_{\Psi_{n,D},S^n,n,\Lambda_{n,D}}\big((g,x),(h,y)\big)
\]
for all $(g,x),(h,y) \in G \times X_{n,D}$ with $h^{-1}g \in S^n$, where $d_{\Psi_{n,D},S^n,n,\Lambda_{n,D}}$ is the metric on $G \times X_{n,D}$ induced from the homotopy coherent $G$-action $\Psi_{n,D}$ on $X_{n,D}$ (Definition \ref{def-metric}).
\end{defn}
\begin{ex}
\begin{enumerate}
\item Every Dress-Farrell-Hsiang group is a DFHJ group
    (Choose  $X_{n,H}$ as a point).
\item If a group $G$ is homotopy transfer reducible over $\mathcal{F}$ then $G$ is a DFHJ group with respect to $\mathcal{F}$ (Choose  $F_n$ as the trivial group).
\end{enumerate}
\end{ex}
\begin{lem} \label{DFHJ-virtual}
Let $G$ be a DFHJ group with respect to the family $\mathcal{F}_G := \{ H < G \mid H \mbox{ satisfies the full Isomorphism Conjecture in A-theory} \}$.
Let $F$ be a finite group. Then the wreath product $G \wr F$ is a DFHJ group with respect to the family $\mathcal{F}_{G \wr F} := \{ H < G \wr F \mid H \mbox{ satisfies satisfies the full Isomorphism Conjecture in A-theory } \}$.
\end{lem}
\Proof The proof follows exactly the same as proof of \cite[Lemma 4.3]{W2} since Dress group is also closed under taking subgroups and quotients.\qed

\begin{thm}\label{DFHJ-iso}
Let $G$ be a DFHJ group with respect to the family   of subgroups $\CF$, then $G$ satisfies the A-theoretic Isomorphism Conjecture with coefficients with respect to the family $\CF$.
\end{thm}

\Proof The proof follows closely to \cite[Proposition 4.9]{W2}. We fix a finite symmetric generating subset $S \subseteq G$ which contains the trivial element $e \in G$. We denote by $d_G$ the word metric with respect to $S \setminus \{e\}$.
Since $G$ is a DFHJ group with respect to $\mathcal{F}$, there exist a natural number $N$ and surjective homomorphisms $\alpha_n \colon G \to F_n$ ($n \in \BN$) with $F_n$ a finite group such that the following condition is satisfied.
For any Dress subgroup $D$ of $F_n$ there exist
\begin{itemize}
\item a compact, contractible, controlled $N$-dominated metric space $X_{n,D}$,
\item a homotopy coherent $G$-action $\Psi_{n,D}$ on $X_{n,D}$,
\item a positive real number $\Lambda_{n,D}$,
\item a simplicial complex $E_{n,D}$ of dimension at most $N$ with a simplicial $\alpha_n^{-1}(D)$-action whose stabilizers belong to $\mathcal{F}$,
\item and an $\alpha_n^{-1}(D)$-equivariant map $f_{n,D} \colon G \times X_{n,D} \to E_{n,D}$
\end{itemize}
such that
\[
n \cdot d_{E_{n,D}}^1\big(f_{n,D}(g,x),f_{n,D}(h,y)\big) \leq d_{\Psi_{n,D},S^n,n,\Lambda_{n,D}}\big((g,x),(h,y)\big)
\]
for all $(g,x),(h,y) \in G \times X_{n,D}$ with $h^{-1}g \in S^n$.
We denote by $\mathcal{D}_n$ the family of Dress subgroups of $F_n$ and $\bar{D}:=\alpha_n^{-1}(D)$. We set

$$S_n := \coprod_{D \in \mathcal{D}_n} G/{\bar{D}}~\times ~G =  \coprod_{D \in \mathcal{D}_n} G~\times_{\bar{D}}~ G$$
 $$Y_n := \coprod_{D \in \mathcal{D}_n} G/{\bar{D}}~\times~ G \times~ X_{n,D} =  \coprod_{D \in \mathcal{D}_n} G~\times_{\bar{D}}~ G \times~ X_{n,D} $$
  $$\Sigma_n :=  \coprod_{D \in \mathcal{D}_n} G~\times_{\bar{D}}~   E_{n,D} \times G$$
Here for example $G~\times_{\bar{D}}~ G$  means $G~\times~ G$ module out the relation $(g,h) = (gr^{-1},rh)$, $r\in \bar{D}$ and we can identify $G~\times_{\bar{D}}~ G$ with $ G/{\bar{D}}~\times ~G$ via sending $(g,g')$ to $(g\bar{D},gg')$ (inverse given by mapping $(gD,g')$ to $(g,g^{-1}g')$).
We will use the quasi-metrics on $ S_n=  \coprod_{D \in \mathcal{D}_n} G~\times_{\bar{D}}~ G$ and $ Y_n =  \coprod_{D \in \mathcal{D}_n} G~\times_{\bar{D}}~ G \times~ X_{n,D} $ given by
\begin{eqnarray*}
d_{ G}\big((g,h)_D,(g',h')_{D'}\big) & := & d_G(gh,g'h'),\\
d_{Y_n}\big((g,h,x)_D,(g',h',x')_{D'}\big) & := & d_{\Psi_{n,D},S^n,n,\Lambda_{n,D}}((gh,x),(g'h',x')),
\end{eqnarray*}
if $D=D'$, $g\bar{D}=g'\bar{D'}$. Otherwise, we set
\[
d_{ S_n}((g,h)_D,(g',h')_{D'}) := \infty \quad \text{and} \quad d_{Y_n}((g,h,x)_D,(g',h',x')_{D'}) := \infty.
\]

If we define the quasi-metric on $\coprod_{D \in \mathcal{D}_n} G/D$ by assigning distance $\infty$ to any points that are not equal. Then $d_{S_n} = d_{\coprod_{D \in \mathcal{D}_n} G/D} +d_G$ and $ d_{Y_n} = d_{\coprod_{D \in \mathcal{D}_n} G/D} +d_{\Psi_{n,D},S^n,n,\Lambda_{n,D}} $.
We define the action of $G$ on $ S_n=  \coprod_{D \in \mathcal{D}_n} G~\times_{\bar{D}}~ G$ via $r(g,h) = (rg,h)$  and the action on  $ Y_n =  \coprod_{D \in \mathcal{D}_n} G~\times_{\bar{D}}~ G \times~ X_{n,D} $ by $r(g,h,x) = (rg,h,x)$ for any $r\in G$.

Similarly,  we define the quasi-metric on $ \Sigma_n$ to be  $ nd^{l^1}_{G\times_{\bar{D}}E_{n,D}} +d_G$ and $G$ acts on $\Sigma_n = \coprod_{D \in \mathcal{D}_n} G~\times_{\bar{D}}~   E_{n,D} \times G$ by $r(g,x,h)=(rg,x,rh)$. We define the map
\[
f_{n} := \coprod_{D \in \mathcal{D}_n} G~\times_{\bar{D}}~ G \times~ X_{n,D} \to \coprod_{D \in \mathcal{D}_n} G~\times_{\bar{D}}~   E_{n,D} \times G ~\colon~ ~ (g,h,x)\to (g,f_{n,D}(h,x),gh)
\]

By Proposition~\ref{obstruction}, it suffices to show the K-theory of the obstruction category
\[
(\mathfrak{R}^G_f(W, \BJ(G,E_{\mathcal{F}}G)),h)
\]
vanishes. This is a direct consequence of the following commuting diagram (abbreviate $E_\mathcal{F} := E_{\mathcal{F}}G$)

  \[
   \begin{tikzpicture}
    \matrix (m) [matrix of math nodes, column sep=4em, row sep=2em, text depth=.5em, text height=1em, ampersand replacement=\&]
    {
    (\mathfrak{R}^G_f(W, \BJ(G,E_{\mathcal{F}})),h)  \& (\mathfrak{R}_{f}^G(W, \BJ(( S_n )_n, E_{\mathcal{F}})), h) \\
    (\mathfrak{R}^G_{f}(W, \BJ((G)_n,E_{\mathcal{F}})),h) \&   (\mathfrak{R}_{fd}^G(W, \BJ(( S_n)_n, E_{\mathcal{F}})), h) \\
      (\mathfrak{R}_{fd}^G(W, \BJ(( G)_n, E_{\mathcal{F}})), h) \& (\mathfrak{R}_{fd}^G(W, \BJ(( S_n)_n, E_{\mathcal{F}})), h^{fin})       \\
    (\mathfrak{R}_{fd}^G(W, \BJ((G)_n, E_{\mathcal{F}})), h^{fin}) \&(\mathfrak{R}_{fd}^G(W, \BJ(( Y_n )_n, E_{\mathcal{F}})), h^{fin})\\
    (\mathfrak{R}_{fd}^G(W, \BJ(( G)_n, E_{\mathcal{F}})), h^{fin}) \& (\mathfrak{R}_{fd}^G(W, \BJ((\Sigma_n)_n, E_{\mathcal{F}})), h^{fin}) \\};
    \path[->]
    (m-1-2) edge node[below]{$P_{j}$}  (m-1-1)
    (m-1-1) edge node[right]{$\Delta$}  (m-2-1)
    (m-2-2) edge node[right]{$i_2$}  (m-3-2)
    (m-4-2) edge node[above]{$P_{ Y_n \rightarrow G}$} (m-4-1)

     (m-1-2) edge node[below]{$~~P^{f}_{ S_n \rightarrow G}$}  (m-2-1)
    (m-3-2) edge node[above]{$~~P_{ S_n \rightarrow G}$}  (m-4-1)
    (m-5-2) edge node[above]{$~~P_{ \Sigma_n \rightarrow G}$}  (m-5-1)
    (m-2-1) edge node[right]{$incl_1$} (m-3-1)
    (m-3-1) edge node[right]{$i_1$} (m-4-1)
    (m-1-2) edge node[right]{$incl_2$} (m-2-2)
     (m-4-1) edge node[right]{$=$}  (m-5-1)
    (m-4-2) edge node[right]{$F$} (m-5-2);
  \draw[dashed, ->]
   (m-1-1) edge [out= 30,in=150] node[below]{$tr_1$}  (m-1-2)
   (m-1-2) to[out= 0,in=0] node[right]{$tr_2$} (m-4-2);
   \end{tikzpicture}
  \]


Here the map $incl_1$, $incl_2$, $i_1$ and $i_2$ are all   inclusions of categories. $\Delta$ is the diagonal map.  The  maps $P_{ \ast \rightarrow G}$ (with or without $f$ decoration) are induced by the obvious  projections. The map $P_{j}$  is defined to be the projection functor which takes the inclusion into the full product category, projects onto the $j$-th component, and then apply the function induced by the projections $\coprod_{D \in \mathcal{D}_j} G/{\bar{D}}~ \times G\times E_{\mathcal{F}}\times [0,\infty[\rightarrow G\times E_{\mathcal{F}}\times [0,\infty[$.

The theorem now follows once we show the following

\begin{enumerate}[label=(\roman*)]

\item In map $incl_1$ and $incl_2$ induces isomorphism after applying K-theory;

\item  The composition $i_1\circ incl_1 \circ \Delta$ induces  injective map  on K-theory;

\item The  K-theory of $(\mathfrak{R}_{fd}^G(W, \BJ((\Sigma_n)_n, E_{\mathcal{F}})), h^{fin})$ vanishes.

\item The map $tr_1$ exists after applying K-theory and the composition $P_{j} \circ tr_1$ induces identity in K-theory  for any $j$, in particular $tr_1$ induces injective maps in K-theory. Moreover, the composition $P^f_{ S_n \rightarrow G}  \circ tr_1$ and $\Delta $ induces the same map in K-theory;

\item The functor $F$, defined as a restriction of
$$
\prod_{n\in\BN}f_n :(\prod_{n\in\BN}\mathfrak{R}_{fd}^G(W, \BJ(( Y_n)_n, E_{\mathcal{F}})), h^{fin} )\rightarrow (\prod_{n\in\BN}\mathfrak{R}_{fd}^G(W, \BJ((  \Sigma_n)_n, E_{\mathcal{F}})), h^{fin} ),
$$
is well defined.
\item The map $tr_2$ exists after applying K-theory. Moreover, the composition  $ P_{ Y_n \rightarrow G} \circ tr_2 $ and $ P_{ S_n \rightarrow G} \circ i_2 \circ  incl_2$ induces the same map in K-theory.

\end{enumerate}
We will prove these statements in positive degrees, for non-positive ones, one has to consider further deloopings of the category, but the proof works in the same way, see \cite[section 5]{UW} and \cite[Section 6.6]{ELPUW} for more details.
We proceed to prove the positive degree case as follows

\begin{enumerate}[label=(\roman*)]
    \item This is stated in \cite[Remark 5.5]{UW};

    \item   This is proved in \cite[Lemma 6.18]{ELPUW} notice that $\Delta$ and $incl_1$ commutes;

     \item  This is proved in \cite[Proposition 6.15(iii)]{ELPUW}, which follows from the ``squeezing Theorem" in A-theory \cite[Theorem 10.1]{UW};

    \item The first part is proved in \cite[Corollary 9.3]{UW}. For the second part, one can argue the same as in \cite[Proposition 9.2(2)]{UW} note that the map $P^f_{ S_n \rightarrow G} $ only changes the controlled conditions on the object.

    \item This is \cite[Proposition 6.15(ii)]{ELPUW};

    \item  The map $tr_2$ can be  defined as in \cite[Section 7]{ELPUW}  with some modifications. We explain now how this should proceed based on the terminology and proof there. Firstly by \cite[Remark 7.3]{ELPUW}, it suffices to define a transfer functor
    $$tr_2^{\alpha,d} :(\mathfrak{R}_{f}^G(W, \BJ((S_n )_n, E_{\mathcal{F}}))_{\alpha,d}, h)  \rightarrow (\mathfrak{R}_{fd}^G(W, \BJ(( Y_n)_n, E_{\mathcal{F}})), h^{fin})$$
    such that the corresponding induced diagram is homotopy commutative on K-theory.
    Now note that $S_n =\coprod_{D \in \mathcal{D}_n} G \times_{\bar{D}} G = \coprod_{D \in \mathcal{D}_n} G/{\bar{D}} \times G$. For each $D \in \mathcal{D}_n$, we have homotopy coherent $G$-action $\Psi_{n,D}$ on a space $X_{n,D}$  satisfies conditions of \cite[Definition 6.2]{ELPUW} by Remark \ref{diff-transfer-reducible}. Thus we have a transfer map  \cite[Lemma 7.12]{ELPUW} $ tr_{2}^{\alpha,d}|_{n,D}$ such that $P_{ G \times X_{n,D} \rightarrow G} \circ tr_{2}^{\alpha,d}|_{n,D}$ and $incl_2$ coincides   after applying K-theory \cite[Proposition 7.24]{ELPUW}

  \[
   \begin{tikzpicture}
    \matrix (m) [matrix of math nodes, column sep=4em, row sep=2em, text depth=.5em, text height=1em, ampersand replacement=\&]
    {
      (\mathfrak{R}_{f}^G(W, \BJ(( g\bar{D}\times G) , E_{\mathcal{F}}))_{\alpha,d}, h)  \& \\
       (\mathfrak{R}_{fd}^G(W, \BJ(( g\bar{D}\times G) , E_{\mathcal{F}})), h)   \& (\mathfrak{R}_{fd}^G(W, \BJ(( g\bar{D}\times G \times X_{n,D}) , E_{\mathcal{F}})), h)  \\};
    \path[->]

    (m-2-2) edge node[above]{$P_{ G \times X_{n,D} \rightarrow G}$} (m-2-1)
    (m-1-1) edge node[right]{$incl_2$} (m-2-1);
  \draw[dashed, ->]
   (m-1-1) edge node [above]{$tr_2^{\alpha,d}|_{n,D}$} (m-2-2);
   \end{tikzpicture}
  \]
  Put all these maps together, we get a transfer map (in fact an exact functor by \cite[Proposition 7.25]{ELPUW} with a little modification)
  $$tr_2^{\alpha,D} :K(\mathfrak{R}_{f}^G(W, \BJ(( \coprod_{D \in \mathcal{D}_n} G/{\bar{D}}~\times~ G \times~ X_{n,D} )_n, E_{\mathcal{F}})), h)_{\alpha,D} $$
  $$ \hspace{50mm}\longrightarrow K(\mathfrak{R}_{fd}^G(W, \BJ((  \coprod_{D \in \mathcal{D}_n} G/{\bar{D}}~\times~ G )_n, E_{\mathcal{F}})), h^{fin}) $$
 The rest of \cite[Section 7.10]{ELPUW} tells us how to get $tr_2$ from $tr_2^{\alpha,D}$.
\end{enumerate}

\section{Dress subgroups}\label{section-dress}
In this section, we embed our group  $ \BZ[w,\frac{1}{w}] \rtimes_w \BZ$ into a slightly bigger group $\mathcal{O}_w\rtimes_w\BZ$ and study Dress subgroups in its finite quotients.
Part of the ideas  in this section are from \cite[Section 4]{FW3} and \cite[Section 5.6]{W2}.

\subsection{The  ring $\mathcal{O}_w$} We review some algebraic number theory background and notions related to the group $ \BZ[w,\frac{1}{w}] \rtimes_w \BZ$ for the reader's convenience, for more details see \cite[Section 5.1]{W2}.

Let $w \in \bar{\BQ}^{\times}$ be a non-zero algebraic number. Let $\mathcal{O}$ be the ring of integers in the algebraic number field $\BQ(w)$, i.e. $\mathcal{O}\subset \BQ(w)$ is the subring consisting of all elements which are integral. Recall that an algebraic number is integral iff it is a root of a monic polynomial with rational integer coefficients. Note that $\mathcal{O}$ is a Dedekind domain, in particular every nonzero prime ideal is maximal.

 For a prime ideal $\mathfrak{p} \subset \mathcal{O}$ we denote by

 $$ \mathcal{O}_\mathfrak{p} := \{ \frac{x}{y} \mid x,y \in \mathcal{O}, y \not \in \mathfrak{p}  \subset \BQ(w) \}$$
the localization of $\mathcal{O}$ at ${\mathfrak{p}}$. Let $V_{\mathfrak{p}} : \BQ(w)^{\times}\rightarrow \BZ$ be the corresponding valuation, i.e., $x \mathcal{O} = \mathfrak{p}^{V_{\mathfrak{p}}(x)} \mathcal{O}_{\mathfrak{p}}$. We extend the valuation $V_\mathfrak{p}$ to $\BQ(w)$ by the convention $V_\mathfrak{p}(0) = \infty$. Note that $\mathcal{O}_{\mathfrak{p}} = \{x \in \BQ(w)\mid V_{\mathfrak{p}}(x) \geq 0 \}$.

A fractional ideal is a finitely generated $\mathcal{O}$-submodule $\mathfrak{a} \neq 0$ of $\BQ(w)$. Every fractional ideal possesses a unique factorization $\mathfrak{a} = \prod_{\mathfrak{p}} \mathfrak{p}^{\nu_{\mathfrak{p}}}$ with $\nu_{\mathfrak{p}} \in \BZ$ and $\nu_{\mathfrak{p}} =0$ for almost all prime ideals $\mathfrak{p}$. We define $M_\alpha$ to be the set of prime ideals appeared in the prime factorization of $\alpha$. For $\mathfrak{a}= x \mathcal{O}$ with $x \in {\BQ(w)}^\times$ we have $\nu_{\mathfrak{p}} = V_{\mathfrak{p}}(x)$. We conclude that

$$M_{x} := M_{x\mathcal{O}} =\{\mathfrak{p} \subset \mathcal{O}  ~prime~ideal\mid V_{\mathfrak{p}}(x)\neq 0 \}$$

is a finite set. In particular, $M_w$ is a finite set.  We define the ring

$$ \mathcal{O}_w  = \{ x \in \BQ(w)\mid V_{\mathfrak{p}}(x)\geq 0 ~for ~all~ prime~ ideals ~\mathfrak{p}\not\in M_w \}$$

Note that $\mathcal{O} \subset \mathcal{O}_{w}$ and $w,w^{-1} \in \mathcal{O}_w$. The group of units in the ring $\mathcal{O}_w $ is given by
$$\mathcal{O}_w^\times =   \{ x \in \BQ(w)\mid V_{\mathfrak{p}}(x)= 0 ~for ~all~ prime~ ideals ~\mathfrak{p}\not\in M_w \}  $$

\subsection{Order of $w$ in a finite quotient}
We  review some results in \cite[Section 5.6]{W2}.

Given any ideal $I$ in $\mathcal{O}$ such that $M_I \cap M_w =\emptyset$, we define a nature number $t_w(I,s)$ for any $s\geq 1$ by the following:
$$t_w(I,s)\BZ = \{ z \in \BZ \mid w^z \equiv 1  \mod I^s\mathcal{O}_w  \} $$

In particular, let $q$ be a prime number and $\mathfrak{q} \not \in M_w$ be a prime ideal in $\mathcal{O}$ which contains $q$, we have

$$t_w(q,s)\BZ = \{ z \in \BZ \mid w^z \equiv 1  \mod q^s\mathcal{O}_w  \} $$
$$t_w(\mathfrak{q},s)\BZ = \{ z \in \BZ \mid w^z \equiv 1  \mod \mathfrak{q}^s\mathcal{O}_w  \} $$

Note that if $q \mathcal{O} = \mathfrak{q}_1^{\nu_1}\mathfrak{q}_2^{\nu_2}\cdots \mathfrak{q}_r^{\nu_r}$, then $M_q = \{\mathfrak{q}_1,\mathfrak{q}_2, \cdots \mathfrak{q}_r\}$ and $t_w(q,s)$ is the least common multiple of $t_w(\mathfrak{q}_i^{\nu_i},s) = t_w(\mathfrak{q}_i,s\nu_i)$ for $i=1,2,\cdots,r$.

\begin{lem}\label{order-q-s}\cite[Lemma 5.29]{W2}
Let $q$ be a prime number and $\mathfrak{q} \not \in M_w$ be a prime ideal in $\mathcal{O}$ which contains $q$, then:
\begin{enumerate}[label=(\roman*)]
\item The ring $\mathcal{O}_w/\mathfrak{q}\mathcal{O}_w$ is a finite field of characteristic $q$ which is isomorphic to $\mathcal{O}/\mathfrak{q}\mathcal{O}$.
\item $t_w(\mathfrak{q},1)$ and $q$ are coprime. In particular, $t_w(\mathfrak{q},1) \neq 0$.
\item For every $s\in \BN$ we have $t_w(\mathfrak{q},s+1) = qt_w(\mathfrak{q},s)$ or $t_w(\mathfrak{q},s+1) = t_w(\mathfrak{q},s)$.
\item Let $(a,b) \in \mathcal{O}_w/\mathfrak{q}^s \mathcal{O}_w \rtimes_w \BZ/t_w(\mathfrak{q},s)$ such that $b \not\in t_w(\mathfrak{q},1) \BZ/ t_w(\mathfrak{q},s)$. Then $(a,b)$ is conjugate to $(0,b)$.
\end{enumerate}

\end{lem}
~
\vspace{5pt}

\subsection{Dress subgroups of finite quotients for $ \mathcal{O}_w \rtimes_w \BZ$} We define our new group $\Gamma_w$ as the semi-product: $ \Gamma_w = \mathcal{O}_w \rtimes_w \BZ$, where the multiplication in $\Gamma_w$ is given by

$$(x_1,y_1)(x_2,y_2) = (x_1 + w^{y_1}x_2,y_1+y_2).$$
for $(x_i,y_i) \in\mathcal{O}_w \rtimes_w \BZ$, $i=1,2$. Since $\BZ[w,\frac{1}{w}] \subset \mathcal{O}$, $\BZ[w,\frac{1}{w}] \rtimes_w \BZ$ is a subgroup of $\Gamma_w$. we want to study its finite quotient.

\begin{lem} \label{primeprogression}
Given any integer $N>0$, there exists  primes $q_1,q_2,q_3>N$  such that
 $M_w$, $M_{q_1}$, $M_{q_2}$,$M_{q_3}$ are pairwise disjoint and $1+w \not \in \mathfrak{q}_i\mathcal{O}_w$ for any $\mathfrak{q}_i \in M_{q_i}(i=1,2,3)$.
\end{lem}

\Proof Given $\mathfrak{p} \in M_w$, then $\mathcal{O}/\mathfrak{p}\mathcal{O}$ is a finite field. Let $M$ be maximal of the characteristic of   these finite fields.  Now choose $q_1,q_2,q_3$ to be different primes that are bigger than $N$ and $M$, they then have the properties in the lemma. Note that there are only finite many prime ideals $\mathfrak{q}$ such that $1+w$ lies in $\mathfrak{q}\mathcal{O}_w$, thus the last condition can also be reached by choosing $q_1,q_2,q_3$ big. \qed

\begin{rem}\label{1+w-invert}
Note that by our choices $1+w$ is invertible in $\mathcal{O}_w/\mathfrak{q}_i^s$ for any $i$ and $s\geq 1$. In particular $1+w$ is invertible in $\mathcal{O}_w/(q_1^sq_2^sq_3^s)$.
\end{rem}

Given any $N>0$, let $q_1,q_2,q_3$ be three primes choosing as in Lemma \ref{primeprogression}.  Recall  $t_w(q_1q_2q_3,s)$ is defined to be the nature number such that
$$t_w(q_1q_2q_3,s)\BZ = \{ z \in \BZ \mid w^z \equiv 1  \mod q_1^sq_2^sq_3^s\mathcal{O}_w  \} $$

Given any $\mathfrak{q}_i \in M_{q_i}$ with $v_{\mathfrak{q_i}}(q_i)=\nu_i\geq 1$, we also consider the nature number $t_w(\mathfrak{q}_1^{\nu_1}\mathfrak{q}_2^{\nu_2}\mathfrak{q}_3^{\nu_3},s)$ defined by

$$t_w(\mathfrak{q}_1^{\nu_1}\mathfrak{q}_2^{\nu_2}\mathfrak{q}_3^{\nu_3},s)\BZ = \{ z \in \BZ \mid w^z \equiv 1  \mod \mathfrak{q}_1^{\nu_1s}\mathfrak{q}_2^{\nu_2s}\mathfrak{q}_3^{\nu_3s}\mathcal{O}_w  \} $$

Note that by Lemma \ref{order-q-s}, $\mathcal{O}_w/\mathfrak{q}_i\mathcal{O}_w$ is a finite field of Characteristic $q_i$, hence $\mathcal{O}_w/\mathfrak{q}_i\mathcal{O}_w$ has order some power of $q_i$. By the Chinese Reminder Theorem,  $\mathcal{O}_w/(\mathfrak{q}_1^{\nu_1s}\mathfrak{q}_2^{\nu_2s}\mathfrak{q}_3^{\nu_3s}) \cong \mathcal{O}_w/\mathfrak{q}_1^{\nu_1s} \oplus  \mathcal{O}_w/\mathfrak{q}_2^{\nu_2s} \oplus \mathcal{O}_w/\mathfrak{q}_3^{\nu_3s} $ and $\mathcal{O}_w / (q_1^sq_2^sq_3^s) \cong \mathcal{O}_w/q_1^s \oplus \mathcal{O}_w/q_2^s \oplus \mathcal{O}_w/q_3^s$. Hence we have $t_w(q_1q_2q_3,s)$ is the least common multiple of $t_w(q_i,s) (i=1,2,3)$ and $t_w(\mathfrak{q}_1^{\nu_1}\mathfrak{q}_2^{\nu_2}\mathfrak{q}_3^{\nu_3},s)$ is the least common multiple of $t_w(\mathfrak{q}_i,\nu_is)(i=1,2,3)$.

Given $n>0$, we define a number $\gamma(s,n)$ depends on $s$ and $n$ by
\begin{equation}\label{gamma-s-n}
\gamma(s,n) = \prod_{i=1,2,3} t_w^{n+1}(\mathfrak{q}_i,s\nu_i)
\end{equation}
Now we define finite quotients of $\mathcal{O}_w \rtimes_w \BZ$ via the following maps.

$$\begin{CD}
\mathcal{O}_w \rtimes_w \BZ  @> {\pi_1} >> \mathcal{O}_w/(\mathfrak{q}_1^{\nu_1s}\mathfrak{q}_2^{\nu_2s}\mathfrak{q}_3^{\nu_3s}) \rtimes_w \BZ \\ @>\pi_2>> \mathcal{O}_w/(\mathfrak{q}_1^{\nu_1s}\mathfrak{q}_2^{\nu_2s}\mathfrak{q}_3^{\nu_3s}) \rtimes_w \BZ/\gamma(s,n) @>Pr>> \BZ/\gamma(s,n)\\
\end{CD}$$

where $\pi_1$ is defined by quotienting out the normal subgroup $ (\mathfrak{q}_1^s\mathfrak{q}_2^s\mathfrak{q}_3^s)  \mathcal{O}_w   \rtimes_w  \{0\}$, $\pi_2$ is defined by quotienting out the normal subgroup $\{ 0 \}  \rtimes_w \gamma(s,n) \BZ$ and $Pr$ is the projection of $\BZ/(q_1^nq_2^nq_3^n) \rtimes_w \BZ/\gamma(s,n)$ to the second factor.

We proceed to study Dress subgroups of $\mathcal{O}_w/(\mathfrak{q}_1^{\nu_1s}\mathfrak{q}_2^{\nu_2s}\mathfrak{q}_3^{\nu_3s}) \rtimes_w \BZ/\gamma(s,n)$.

\begin{lem} \label{prelemma-dress-subgroup-gel}
Given a positive integer $n$, there exist primes $q_1,q_2,q_3 > n$ and a positive integer $s >n$ such that for any $\mathfrak{q}_i \in M_{q_i}(i=1,2,3)$ one of the following is true for each Dress subgroup $D$ of $\mathcal{O}_w/(\mathfrak{q}_1^{\nu_1s}\mathfrak{q}_2^{\nu_2s}\mathfrak{q}_3^{\nu_3s}) \rtimes_w \BZ/\gamma(s,n)$:
\begin{enumerate}[label=(\roman*)]
    \item \label{prelemma1} The index $[\BZ/\gamma(s,n), Pr(D)] \geq n$;
    \item \label{prelemma2} $D \cap Q_i = \{0\}$ for some $i =1,2 ~or~ 3$ where $Q_i$ denotes the $q_i$-Sylow subgroup of $\mathcal{O}_w/(\mathfrak{q}_1^{\nu_1s}\mathfrak{q}_2^{\nu_2s}\mathfrak{q}_3^{\nu_3s})\rtimes_w \{0\}$.
\end{enumerate}
\end{lem}

\Proof Let $q_1,q_2,q_3$ be the three primes determined by Lemma \ref{primeprogression} with $N=n$. We now choose $s$ big enough, so that $w^z-1 \not \in \mathfrak{q}_i^s$ for any $1 \leq z \leq t_w(\mathfrak{q}_i,1)$ for $i=1,2,3$. This implies $\frac{t_w(\mathfrak{q}_i,s)}{t_w(\mathfrak{q}_i,1)} >1$. Since $\frac{t_w(\mathfrak{q}_i,s)}{t_w(\mathfrak{q}_i,1)}$ is a power of $q_i$ by Lemma \ref{order-q-s} (iii) and $q_i>n$, we have $\frac{t_w(\mathfrak{q}_i,s)}{t_w(\mathfrak{q}_i,1)} >n$. Applying Lemma \ref{order-q-s} (iii) again, we have $\frac{t_w(\mathfrak{q}_i,\nu_is)}{t_w(\mathfrak{q}_i,1)} >n$ since $\nu_i\geq 1$ for $i=1,2,3$. In particular $q_i$ appears in the prime factorization of $t_w(\mathfrak{q}_1^{\nu_1}\mathfrak{q}_2^{\nu_2}\mathfrak{q}_3^{\nu_3},s)$.

Now let $D$ be a Dress subgroup with $D_1 \lhd D_2 \lhd D$ such that $D_1$ is a
$p_1$-group, $D_2/D_1$ is cyclic and $D/D_2$ is a $p_2$-group. By Lemma \ref{dressdef}, we can assume $(p_1,|D_2/D_1|) = (p_2,|D_2/D_1|) =1$. Since we know $\gamma(s,n)$'s prime decomposition contains all the prime factor $q_i$ and $q_i>n$, without loss of generality, we can assume $p_1 \neq q_3$ and $p_2 \neq q_3$. In fact, we can assume
$p_1= q_1, p_2=q_2$ otherwise property (i) in the Lemma holds.

Assume now $D_2/D_1 \cong C$, then $D_2 \cong D_1 \rtimes C$, where $C$ is a finite cyclic group. Notice that
$$D \cap Q_3 = C \cap Q_3$$
Now let the generator of $C$ to be $(a,b)\in  \mathcal{O}_w/(\mathfrak{q}_1^{\nu_1s}\mathfrak{q}_2^{\nu_2s}\mathfrak{q}_3^{\nu_3s}) \rtimes_w \BZ/\gamma(s,n)$. Then

$$b \in ( \frac{t_w^{n+1}(\mathfrak{q}_1,\nu_1s)}{t_w^{n+1}(\mathfrak{q}_1,1)} \cdot \frac{t_w^{n+1}(\mathfrak{q}_2,\nu_2s)}{t_w^{n+1}(\mathfrak{q}_2,1)} )  \BZ/  \gamma(s,n)$$
and the order of $b$ divides $t_w^{n+1}(\mathfrak{q}_1,1) {t_w^{n+1}(\mathfrak{q}_2,1)}  {t_w^{n+1}(\mathfrak{q}_3,\nu_3s)} $. Let the order of $b \in \BZ/\gamma(s,n)$ to be $k$, then $D\cap Q_3$ is a cyclic group generated by $(a,b)^k$. If we assume property (i) in the Lemma does not hold,  then

\begin{equation}\label{key-ineq}
(t_w^{n+1}(\mathfrak{q}_1,1) {t_w^{n+1}(\mathfrak{q}_2,1)}  {t_w^{n+1}(\mathfrak{q}_3,\nu_3s)} ) :k <n
\end{equation}

Note that this implies the order of $b$ must divides $t_w(\mathfrak{q}_3,1)$. In fact, if this is not true, then there exist a prime $p$ such that $p^l$ divides $t_w(\mathfrak{q}_3,1)$, but $p^l$ does not divide $k$ for some $l\geq 1$. This implies $p^{(n+1)l}$ divides $t_w^{n+1}(\mathfrak{q}_3,1)$, but $p^l$ does not divide $k$, we have
$$t_w^{n+1}(\mathfrak{q}_3,1):k \geq p^{nl}>n$$
which is a contradiction to the inequality \ref{key-ineq}. Since $q_3>n$, we have $t_w(\mathfrak{q}_3,\nu_3s)$ also divides $k$ otherwise property (i) in the lemma holds.

Now we can assume $b$'s order $k = k't_w(\mathfrak{q}_3,\nu_3s)$, then we have the following calculation in our group $\mathcal{O}_w/(\mathfrak{q}_3^s\mathfrak{q}_2^s\mathfrak{q}_3^s)\rtimes_w \BZ/\gamma(s,n)$:

$$(a,b)^k = ((a,b)^{t_w(\mathfrak{q}_3,\nu_3s)})^{k'}= ((1+w+\cdots w^{t_w(\mathfrak{q}_3,\nu_3s) -1})b ,b)^{k'}$$

Now by our choices of $q_i$ and Remark \ref{1+w-invert}, we have $1+w$ is invertible in $\mathcal{O}_w/(\mathfrak{q}_1^{\nu_1s}\mathfrak{q}_2^{\nu_2s}\mathfrak{q}_3^{\nu_3s})$, thus

$$(a,b)^k =((1+w)^{-1} (1-w^{t_w(\mathfrak{q}_3,\nu_3s)}) b ,b)^{k'}$$

Now $w^{t_w(\mathfrak{q}_3,\nu_3s)} \equiv 1 \mod \mathfrak{q}_3^{\nu_3s}$ in $\mathcal{O}_w$, we have $(1+w)^{-1}(1-w^{t_w(\mathfrak{q}_3,\nu_3s)}) b \in \mathfrak{q}_3^{\nu_3s} \mathcal{O}_w$ and $(a,b)^k \cap Q_3 = \{0\}$.

\qed

Similar to  Lemma \ref{order-q-s} (iv), we have the following

\begin{lem}\label{cong-gel}
Let $q_1,q_2,q_3$ be primes such that $M_w$, $M_{q_1}$, $M_{q_2}$,$M_{q_3}$ are pairwise disjoint and $(a,b) \in \mathcal{O}_w/(\mathfrak{q}_1^{\nu_1s}\mathfrak{q}_2^{\nu_2s}\mathfrak{q}_3^{\nu_3s}) \rtimes_w \BZ/\gamma(s,n)$. If $b \not\in t_w(\mathfrak{q}_i,1) \BZ/ \gamma(s,n)$ for $i=1,2,3$, then $(a,b)$ is conjugate to $(0,b)$.
\end{lem}

\Proof Note first that $t_w(\mathfrak{q},\nu_is)$ divides $\gamma(s,n)$. Now if  $b \not\in t_w(\mathfrak{q}_i,1) \BZ/ \gamma(s,n)$, by the proof of \cite[Lemma 5.29(5)]{W2}, we have $w^b-1$  is a unit in $\mathcal{O}_w/\mathfrak{q}_i^{\nu_is}$ for $i =1,2,3$. Hence $w^b-1$ is a unit in $\mathcal{O}_w/(\mathfrak{q}_1^{\nu_1s}\mathfrak{q}_2^{\nu_2s}\mathfrak{q}_3^{\nu_3s})$. Now set $x:= (w^b-1)^{-1} a \in \mathcal{O}_w/(\mathfrak{q}_1^{\nu_1s}\mathfrak{q}_2^{\nu_2s}\mathfrak{q}_3^{\nu_3s})$. Then $(x,0)(a,b)(x,0)^{-1} = (0,b)$.\qed

We continue to refine Lemma \ref{prelemma-dress-subgroup-gel} to the following using Lemma \ref{cong-gel}.

\begin{prop} \label{prop-dress-subgroup-gel}
Given a positive integer $n$, there exists primes $q_1,q_2,q_3 > n$ and a positive integer $s >n$ such that for any $\mathfrak{q}_i \in M_{q_i}$ one of the following is true for each Dress subgroup $D$ of $\mathcal{O}_w/(\mathfrak{q}_1^s\mathfrak{q}_2^s\mathfrak{q}_3^s)\rtimes_w \BZ/\gamma(s,n)$:
\begin{enumerate}[label=(\roman*)]
    \item \label{thm1} The index $[\BZ/\gamma(s,n), Pr(D)] \geq n$;
    \item \label{thm2} there exists $g\in \mathcal{O}_w/(\mathfrak{q}_1^{\nu_1s}\mathfrak{q}_2^{\nu_2s}\mathfrak{q}_3^{\nu_3s})\rtimes_w \BZ/\gamma(s,n)$ such that  $gDg^{-1} \subset \mathfrak{q}_i^{\nu_i s} \mathcal{O}_w/(\mathfrak{q}_1^{\nu_1s}\mathfrak{q}_2^{\nu_2s}\mathfrak{q}_3^{\nu_3s})\rtimes_w \BZ/\gamma(s,n)$ for some $i =1,2 ~or~ 3$.
\end{enumerate}
\end{prop}

\Proof By Lemma \ref{prelemma-dress-subgroup-gel}, we can choose primes $q_1,q_2,q_3$ and $s$ big such that, if the index $[\BZ/\gamma(s,n), Pr(D)] < n$ then $D \cap Q_i \subset \mathfrak{q}_i^{\nu_i s} Q_i$ for some $i =1,2 ~or~ 3$ where $Q_i$ denotes the $q_i$-Sylow subgroup of $\mathcal{O}_w/(\mathfrak{q}_1^{\nu_1s}\mathfrak{q}_2^{\nu_2s}\mathfrak{q}_3^{\nu_3s}) \rtimes_w \{0\}$. We can further assume $q_i$ big such that $q_i$ lies in none of the prime ideals which appear in the prime factorization of the fractional ideal $(w^l-1)\mathcal{O}$ with $1 \leq l\leq n$ for $i=1,2,3$. This guarantees that $t_w(\mathfrak{q}_i,1)>n$ for $\mathfrak{q}_i \in M_{q_i}$.

Now Consider the projection $pr$ restricted to $D$, then $pr(D)$ is a cyclic subgroup of $\BZ/\gamma(s,n)$ and the kernel $K = D \cap \mathcal{O}_w/(\mathfrak{q}_1^{\nu_1s}\mathfrak{q}_2^{\nu_2s}\mathfrak{q}_3^{\nu_3s}) \rtimes_w\{0\}$. We have the following short exact sequence:

\begin{equation} \label{ses-split}
1 \rightarrow K \rightarrow D \rightarrow Pr(D) \rightarrow 1.
\end{equation}

Choose any $(a,b) \in D$ such that $b$ generates $Pr(D)$. Assume the Dress subgroup $D$ does not satisfy property (i) of our Proposition, then

\textbf{Claim:} The subgroup generated by $(a,b)$ in $D$ is isomorphic to $Pr(D)$ via $Pr$. In particular, $D$ can be generated by elements in $K$ and $(a,b)$.

 Let the order of $b$ to be $k$. If property (i) does not hold,  then $\frac{\gamma(s,n)}{k} \leq n$.  Similar to the arguments in the proof of  Lemma \ref{prelemma-dress-subgroup-gel} (after  inequality \ref{key-ineq}), we have  $t_w(\mathfrak{q}_i,\nu_is)$ divides $k$ for  any $i=1,2,3$. Since $1-w$ is invertible in $\mathcal{O}_w/\mathfrak{q}_i^{\nu_is}$, if we let $k=t_w(\mathfrak{q}_i,\nu_is)k_i$, then

$$ (a,b)^k =((a,b)^{t_w(\mathfrak{q}_i,\nu_is)})^{k_i}=((1-w)^{-1}(1-w^{t_w(\mathfrak{q}_i,\nu_is)}),t_w(\mathfrak{q}_i, \nu_is)b)^{k_i} $$

Now $w^{t_w(\mathfrak{q}_i,\nu_is)} = 1 \in \mathcal{O}_w/\mathfrak{q}_i^{\nu_is}$ for $i=1,2,3$. We have
$$(a,b)^k = (0,0)\in \mathcal{O}_w/(\mathfrak{q}_1^{\nu_1s}\mathfrak{q}_2^{\nu_2s}\mathfrak{q}_3^{\nu_3s})\rtimes_w \BZ/\gamma(s,n).$$
This proves the claim.  Now $b \not \in t_w(\mathfrak{q}_i,1)\BZ/\gamma(s,n)$ for any $i$ otherwise the index $[\BZ/\gamma(s,n), Pr(D)] \geq n$ since by our choices $t_w(\mathfrak{q}_1,1)>n$. This means we can conjugate $(a,b)$ to $(0,b)$ by Lemma \ref{cong-gel} via  $g = (w^b-1)^{-1} a,0)$. Since $gKg^{-1} \in \mathcal{O}_w/(\mathfrak{q}_1^{\nu_1s}\mathfrak{q}_2^{\nu_2s}\mathfrak{q}_3^{\nu_3s}) \rtimes_w\{0\}$, we have $gDg^{-1}$ is generated by elements in $\mathcal{O}_w/(\mathfrak{q}_1^{\nu_1s}\mathfrak{q}_2^{\nu_2s}\mathfrak{q}_3^{\nu_3s}) \rtimes_w \{0\}$ and $(0,b)$. Now $gDg^{-1} \cap Q_i  = gKg^{-1}\cap Q_i =\{0\}$ for some $i =1,2 ~or~ 3$ where $Q_i$ denotes the $q_i$-Sylow subgroup of $\mathcal{O}_w/(\mathfrak{q}_1^{\nu_1s}\mathfrak{q}_2^{\nu_2s}\mathfrak{q}_3^{\nu_3s}) \rtimes_w \{0\}$, hence $gDg^{-1} \cap \mathcal{O}_w/(\mathfrak{q}_1^{\nu_1s}\mathfrak{q}_2^{\nu_2s}\mathfrak{q}_3^{\nu_3s}) \rtimes_w\{0\} \subset \mathfrak{q}_i^{\nu_i s} \mathcal{O}_w/(\mathfrak{q}_1^{\nu_1s}\mathfrak{q}_2^{\nu_2s}\mathfrak{q}_3^{\nu_3s}) \rtimes_w\{0\} $ and we have now Property (ii) in the theorem holds.

\qed

\begin{rem} \label{assumption-gel}
We  summarize choices we need to make for $q_1,q_2,q_3$ and $s$ such that Proposition \ref{prop-dress-subgroup-gel} holds:
\begin{enumerate}
\item $q_1,q_2,q_3$ are different primes bigger than $n$;
\item $M_w$, $M_{q_1}$, $M_{q_2}$ and $M_{q_3}$ are pairwise disjoint;
\item $1+w \not \in \mathfrak{q}_i\mathcal{O}_w$ for any $\mathfrak{q}_i \in M_{q_i}(i=1,2,3)$.
\item $q_i$ lies in none of the prime ideals which appear in the prime factorization of the fractional ideal $(w^l-1)\mathcal{O}$ with $1 \leq l\leq n$ for $i=1,2,3$;
\item \label{last-c}   Choose $s$ big such that $w^z-1 \not \in \mathfrak{q}_i^s$ for any $1 \leq z \leq t_w(\mathfrak{q}_i,1)$ for $i=1,2,3$;
\end{enumerate}
\end{rem}
\begin{rem}\label{choose-g}
For any $(a,b)\in D$ such that $b$ generates $Pr(D)$, the proof shows the $g$ in (ii) can be chosen to be $((w^b-1)^{-1} a,0)$.
\end{rem}

We proceed to prove our main theorem in this section. Similar to $\gamma(s,n)$, we define a number $\bar\gamma(s,n)$ depends on $s$ and $n$ by
\begin{equation}\label{bargamma-s-n}
\bar\gamma(s,n) = \prod_{i=1,2,3} t_w^{n+1}(q_i,s)
\end{equation}

And we define another finite quotient of $\mathcal{O}_w \rtimes_w \BZ$ via the following maps.

$$\begin{CD}\label{bar-gamma-s-n}
\mathcal{O}_w \rtimes_w \BZ  @> {\bar\pi_1} >> \mathcal{O}_w/(q_1^sq_2^sq_3^s) \rtimes_w \BZ @>\bar\pi_2>> \mathcal{O}_w/(q_1^sq_2^sq_3^s) \rtimes_w \BZ/\bar\gamma(s,n) @>\bar{Pr}>> \BZ/\bar\gamma(s,n)\\
\end{CD}$$

where $\bar\pi_1$ is defined by quotienting out the normal subgroup $ (q_1^sq_2^sq_3^s)  \mathcal{O}_w   \rtimes_w  \{0\}$, $\bar\pi_2$ is defined by quotienting out the normal subgroup $\{ 0 \}  \rtimes_w \bar\gamma(s,n) \BZ$ and $\bar{Pr}$ is the projection of $\BZ/(q_1^nq_2^nq_3^n) \rtimes_w \BZ/\bar\gamma(s,n)$ to the second factor. 

The following theorem is parallel to \cite[Corollary 5.31]{W2}.

\begin{thm} \label{thm-dress-subgroup-gel}
Given a positive integer $n$, there exists primes $q_1,q_2,q_3 > n$ and a positive integer $s >n$ such that
 one of the following is true for each Dress subgroup $D$ of $\mathcal{O}_w/(q_1^s q_2^s q_3^s)\rtimes_w \BZ/\bar \gamma(s,n)$:
\begin{enumerate}[label=(\roman*)]
    \item \label{thm1} The index $[\BZ/\bar \gamma(s,n), \bar{Pr}(D)] \geq n$;
    \item \label{thm2} there exists $g\in \mathcal{O}_w/(q_1^sq_2^sq_3^s)\rtimes_w \BZ/\bar \gamma(s,n)$ such that  $gDg^{-1} \subset q_i^s \mathcal{O}_w/(q_1^sq_2^sq_3^s)\rtimes_w \BZ/\bar\gamma(s,n)$ for some $i =1,2 ~or~ 3$.
\end{enumerate}
\end{thm}

\Proof
We choose $q_1,q_2,q_3$ under the assumptions of Remark \ref{assumption-gel} with condition (\ref{last-c}) strengthened to the following: Choose $s$ big such that $w^z-1 \not \in \mathfrak{q}_i^s$ for $\mathfrak{q}_i\in M_{q_i}$ and any $1 \leq z \leq t_w(\mathfrak{q}_i,1)$ for $i=1,2,3$. Then Proposition \ref{prop-dress-subgroup-gel} holds for any such $q_1,q_2,q_3$ and any $\mathfrak{q}_i\in M_{q_i}(i=1,2,3)$. Now let $q_i$'s prime ideal factorization be the following:

$$q_i \mathcal{O}= \mathfrak{q}_{i1}^{\nu_{i1}} \mathfrak{q}_{i2}^{\nu_{i2}} \cdots \mathfrak{q}_{ir_i}^{\nu_{ir_i}}$$

By the Chinese Remainder Theorem and the fact that $M_w$, $M_{q_1}$, $M_{q_2}$ and $M_{q_3}$ are pairwise disjoint, we have

$$\mathcal{O}_w/(q_1^sq_2^sq_3^s) \cong \bigoplus_{1\leq i\leq 3} \mathcal{O}_w/q_i^s \cong  \bigoplus_{1\leq i\leq 3,1\leq j\leq r_i} \mathcal{O}_w/ \mathfrak{q}_{ij}^{\nu_{ij}} $$

Now let $D$ be a Dress subgroup with $D_1 \lhd D_2 \lhd D$ such that $D_1$ is a
$p_1$-group, $D_2/D_1$ is cyclic and $D/D_2$ is a $p_2$-group. We assume property (i) in the theorem does not hold.  Similar to the proof of Lemma \ref{prelemma-dress-subgroup-gel}, without loss generality,  we can assume
$p_1= q_1, p_2=q_2$.

Now Consider the projection $\bar{pr}$ restricted to $D$, then $\bar{pr}(D)$ is a cyclic subgroup of $\BZ/\bar\gamma(s,n)$ and the kernel $\bar{K} = D \cap \mathcal{O}_w/(q_1^sq_2^sq_3^s) \rtimes_w\{0\}$. We have the following short exact sequence:

\begin{equation} \label{ses-split-j}
1 \rightarrow \bar{K} \rightarrow D \rightarrow \bar{Pr}(D) \rightarrow 1.
\end{equation}

Now choose $(a,b)\in D$ such that $b$ generates $\bar{Pr}(D)$. Note that $w^b-1 \not \in \mathfrak{q}_{i} \mathcal{O}_w$ for any $\mathfrak{q}_i \in M_{q_i}$, otherwise $b$ divides $t(\mathfrak{q}_{ij},1)$ and we have property (i) holds. This means $w^b-1$ is invertible in $\mathcal{O}_w/(q_1^sq_2^sq_3^s)$. We claim now, if we let $g = ((w^{b}-1)^{-1} a,0)$, then $gDg^{-1} \subset q_3^s \mathcal{O}_w/(q_1^sq_2^sq_3^s)\rtimes_w \BZ/\bar\gamma(s,n)$.

In fact, fix any two prime ideals $\mathfrak{q}_1 \in M_{q_1}$, $\mathfrak{q}_2 \in M_{q_2}$. For each $1\leq j\leq r_3$, we define projections via the following maps

$$\begin{CD}
\mathcal{O}_w/(q_1^sq_2^sq_3^s) \rtimes \BZ/\bar{\gamma}(s,n)  @> {\rho_{1j}} >> \mathcal{O}_w/(\mathfrak{q}_1^{\nu_1s}\mathfrak{q}_2^{\nu_2 s}\mathfrak{q}_{3j}^{\nu_{3j}s}) \rtimes_w \bar{\gamma}(s,n) \\
 @>\rho_{2j}>> \mathcal{O}_w/(\mathfrak{q}_1^{\nu_1s}\mathfrak{q}_2^{\nu_2 s}\mathfrak{q}_{3j}^{\nu_{3j}s}) \rtimes_w \BZ/\gamma_j(s,n) @> {\bar{Pr}_j} >> \BZ/\gamma_j(s,n)\\
\end{CD}$$
where $\gamma_j(s,n)$ is the $\gamma(s,n)$ corresponding to the prime ideals $\mathfrak{q}_1,\mathfrak{q}_2$ and $\mathfrak{q}_{3j}$ defined by equation \ref{gamma-s-n}. We set $\rho_j = \rho_{2j}\circ \rho_{1j}$ and  $\rho_j(a,b)=(a_j,b_j)\in \mathcal{O}_w/(\mathfrak{q}_1^{\nu_1s}\mathfrak{q}_2^{\nu_2 s}\mathfrak{q}_{ij}^{\nu_{3j}s}) \rtimes_w \BZ/\gamma_j(s,n)$.

Now a quotient of a Dress group is still a Dress group, thus  $\rho_j(D)$ is a Dress subgroup of $\mathcal{O}_w/(\mathfrak{q}_1^{\nu_1s}\mathfrak{q}_2^{\nu_2 s}\mathfrak{q}_{3j}^{\nu_{3j}s}) \rtimes_w \BZ/\gamma_j(s,n)$. Consider the projection $\bar{Pr}_j$ restricted to $\rho_j(D)$, then  since we assume property (i) in the theorem does not hold,  we have the index $[\BZ/\gamma_j(s,n),\rho_j(D)] \leq n$. By the claim in the proof of Proposition \ref{prop-dress-subgroup-gel}, we have $\bar{pr}_j(\rho_j(D))$ is a cyclic subgroup of $\BZ/\gamma_j(s,n)$ generated by $b_j$ and $(a_j,b_j)$ generate a cyclic subgroup in $\rho_j(D)$ isomorphic to $\bar{pr}_j(\rho_j(D))$ via  $\bar{pr}_j$.

Now by Proposition \ref{prop-dress-subgroup-gel} and Remark \ref{choose-g}, we have

$$((w^{b_j}-1)^{-1} a_j,0)\rho_j(D)((w^{b_j}-1)^{-1} a_j,0)^{-1} \in \mathfrak{q}_3^{\nu_{3j}s} \mathcal{O}_w/(\mathfrak{q}_1^{\nu_1s}\mathfrak{q}_2^{\nu_2 s}\mathfrak{q}_{3j}^{\nu_{3j}s}) \rtimes_w \BZ/\gamma_j(s,n)$$

And since $\rho_j(((w^{b}-1)^{-1} a,0))=((w^{b_j}-1)^{-1} a_j,0)$,  we have for any $1\leq j\leq r_3$

$$((w^{b}-1)^{-1} a,0)D((w^{b}-1)^{-1} a,0)^{-1} \subset \mathfrak{q}_{3j}^{\nu_{3j}s} \mathcal{O}_w/(q_1^sq_2^sq_3^s)\rtimes_w \BZ/\bar\gamma(s,n)$$

Therefore,
$$((w^{b}-1)^{-1} a,0)D((w^{b}-1)^{-1} a,0)^{-1} \subset q_3^s \mathcal{O}_w/(q_1^sq_2^sq_3^s)\rtimes_w \BZ/\bar\gamma(s,n)$$

\qed

\begin{rem}\label{lift-conjugation}
Since the map $\bar\pi:\mathcal{O}_w \rtimes_w \BZ \rightarrow \mathcal{O}_w/(q_1^sq_2^sq_3^s) \rtimes_w \BZ/ \bar{\gamma}(s,n) $ is surjective, our $g$ can be lifted to an element in  $\mathcal{O}_w \rtimes_w \BZ$.
\end{rem}

\section{Proof of the main theorem}\label{proof}

We prove our main theorem in this section.

\begin{lem}\label{reduction}
The full Isomorphism Conjecture in A-theroy holds for solvable groups iff it holds for $G_w=\BZ[w,\frac{1}{w}]\rtimes_w \BZ$ where $w$ is any non-zero algebraic number.
\end{lem}

\Proof The proof follows exactly the same as \cite[Proposition 3.3]{W2}. In fact, the proof there only uses inheritance properties of the Isomorphism Conjecture in K- and L-theory and the fact that the conjecture holds for hyerbolic and CAT(0) groups. By \cite[Theorem 1.1]{ELPUW}, these  inheritance properties and results also hold for the Isomorphism Conjecure in A-theory. \qed

\begin{prop}\label{nearly-crystal}
Let $w \in \bar{\BQ}^{\times}$ be a non-zero algebraic number of infinite order. Then the group $G_w$ is a DFHJ group with respect to the family of subgroups that satisfy the full Isomorphism Conjecture in A-theory.
\end{prop}

\Proof The proof again follows almost exactly the same as in \cite[Proposition 5.33]{W2}. Notice first that the only difference between the definition of Dress-Farrell-Hsiang-Jones group (Definition \ref{def-DFHJ}) and Farrell-Hsiang-Jones group \cite[Definition 4.1]{W2} is that one has to replace the hyper-elementary subgroup to Dress-subgroup. This issue is addressed already in Theorem \ref{thm-dress-subgroup-gel}. We sketch now how should the proof actually go based on the proof of \cite[Proposition 5.33]{W2}.

We choose $N$ to be the number appearing in \cite[Proposition 5.26]{W2}. Given $n>0$, let $\bar{n} = 4nm_2$ and $s> \max\{1,-\log_2 B\}$, where $m_2$ is defined in the beginning of  proof of \cite[Proposition 5.33]{W2} and $B$ is defined in \cite[Proposition 5.26]{W2}.

Now choose $q_1,q_2,q_3$ and $s$ as in Theorem \ref{thm-dress-subgroup-gel} for $n=\bar{n}$ and define $\alpha_n$ as the composition of the following maps
$$ \BZ[w,\frac{1}{w}]\rtimes_w \BZ \hookrightarrow \mathcal{O}_w\rtimes_w \BZ \rightarrow \mathcal{O}_w/(q_1^s q_2^s q_3^s)\rtimes_w \BZ/\bar \gamma(s,\bar{n})$$
in particular $F_n =   \mathcal{O}_w/(q_1^s q_2^s q_3^s)\rtimes_w \BZ/\bar \gamma(s,\bar{n})$, where $\bar\gamma(s,\bar{n})$ was defined by Equation \ref{bargamma-s-n} and quotient map are defined right below that. By Theorem \ref{thm-dress-subgroup-gel}, we have for every Dress subgroup $D$ in $F_n$
\begin{itemize}
    \item Case (1), there exists $g\in \mathcal{O}_w/(q_1^sq_2^sq_3^s)\rtimes_w \BZ/\bar \gamma(s,n)$ such that  $gDg^{-1} \subset q_i^s \mathcal{O}_w/(q_1^sq_2^sq_3^s)\rtimes_w \BZ/\bar\gamma(s,n)$ for some $i =1,2 ~or~ 3$, note that $q_i^{-s} \leq 2^{-s} \leq B$ and $g$ can be lifted to an element in $\mathcal{O}_w\rtimes_w \BZ$ (Remark \ref{lift-conjugation}).
     \item  Case (2), The index $[\BZ/\bar \gamma(s,n), \bar{Pr}(D)] \geq \bar{n} = 2nm_2$;
\end{itemize}

The rest of the proof follows as in \cite[Proposition 5.33]{W2} without change. In particular, for case (1), we choose the corrsponding space $X_{n,D}$, number $\Lambda_{n,D}$ and simplicial complex $E_{n,D}$ via \cite[Proposition 5.26]{W2}. For case (2), we choose $X_{n,D}$ to be a point and $E_{n,D}$ to be the real line with elements $l\BZ$ as vertices, where $l = [\BZ/\bar \gamma(s,n), \bar{Pr}(D)] $. \qed

 Now when $w$ is a non-zero algebraic number of infinite order, by  Lemma \ref{DFHJ-virtual},  any finite wreath product of $G_w$  is also a DFHJ group with respect to the family of subgroups that already satisfy the full Isomorphism Conjecture. By Theorem \ref{DFHJ-iso} and the transitivity principle \cite[Proposition 11.2]{UW}, we have $G_w$  satisfies the full Isomorphism Conjecture  in A-theory. On the other hand, if $w$ is a non-zero algebraic number of finite order, then $G_w$ is virtually abelian and hence satisfies the full Isomorphism Conjecture by \cite[Theorem 1.1]{UW}. Thus by Lemma \ref{reduction}, we have the full Isomorphism Conjecture in A-theory  holds for all solvable groups.

\end{document}